\documentclass{article}

\usepackage{arxiv}

\usepackage[utf8]{inputenc} 
\usepackage[T1]{fontenc}    
\usepackage{hyperref}       
\usepackage{url}            
\usepackage{booktabs}       
\usepackage{amsfonts}       
\usepackage{nicefrac}       
\usepackage{microtype}      
\usepackage{lipsum}
\usepackage{graphicx}  
\usepackage{amsmath}     
\usepackage{amssymb}
\usepackage{bbold}
\usepackage[dvips]{epsfig}
\title{Comments on the Bellman functional for linear time-delay systems}

\author{
 Jorge-Manuel Ortega-Mart\'inez\\
  Department of Automatic Control\\ 
  CINVESTAV-IPN\\
  Ciudad de M\'exico, M\'exico \\
  \texttt{jortega@ctrl.cinvestav.mx} \\
   \And
 Omar-Jacobo Santos-S\'anchez \\
  CITIS-ICBI\\
  Autonomous Hidalgo State University\\
  Carretera Pachuca-Tulancingo,\\
  Mineral de la Reforma, Hidalgo, M\'exico \\
  \texttt{omarj@uaeh.edu.mx} \\
  \And
 Sabine Mondi\'e\\
  Department of Automatic Control\\ 
  CINVESTAV-IPN\\
  Ciudad de M\'exico, M\'exico \\
  \texttt{smondie@ctrl.cinvestav.mx} \\
}
\newtheorem{thm}{Theorem}
\newtheorem{prop}{Proposition}

\begin{document}
\maketitle

\begin{abstract}
In this note, we present some complementary results on the infinite horizon optimal control for linear time-delay systems. We formally establish some properties of the matrices arising in the Bellman functional, and we prove that no concentrated delay term is present in the optimal control law.
\end{abstract}


\section{Introduction}
The optimal control problem for linear time-delay systems  was studied for the first time by Krasovskii \cite{Krasovskii_1962,Krasovskii1963} in the framework of Dynamic Programming. There, sufficient stability conditions were given and the general form of the Bellman functional was suggested. This functional was the starting point for the explicit characterization of the optimal control  was presented a few years later in Ross \cite{Ross_1969}. A numerical example  was exposed in \cite{Ross71}. However, neither this form for the Bellman functional, nor some of its properties were formally justified. In our recent work, based on the strategy introduced in \cite{Santos_2009} relying on the fundamental and Lyapunov matrix of the delay system, \cite{kharitonov2012}, we have been able to give in our submitted contribution \cite{Automatica2020}, formal arguments for the choice of the form of the Bellman functional: In the careful review process of the contribution, a number of pertinent questions and  comments were made. Some of them seemed important to us but could not be answered thoroughly in the paper due to length limitations, so we answer them in this supplementary material.
\\

The note is organized as follows: We briefly recapitulate the known results on the optimal control of time-delay systems in section 2. In section 3 we give the long proof of some properties of the matrices appearing in the Bellman functional, and in section 4 we justify the absence of a concentrated delay term in the optimal control law. We end the note with some short concluding remarks.

We denote the space of $\mathbb{R}^{n}$-valued piecewise-continuous functions on $[-h, 0]$ by $\mathrm{PC}([-h, 0],\mathbb{R}^{n})$. For a given initial function $\varphi(\theta)$,  $x_{t}(\varphi)$ denotes the state of the delay system $\{ x(t+\theta , \varphi), \theta \in [-h,0]\}$, with delay $h>0$; when the initial condition is not crucial, the argument $\varphi$ is omitted. The Euclidian norm for vectors is represented by $\parallel \cdot \parallel$. The set of piecewise continuous functions is equipped with the norm $\parallel \varphi \parallel_{h}=\sup_{\theta\in [-h,0]}\parallel \varphi(\theta) \parallel$. The notations $Q>0$ indicates that matrix $Q$ is positive definite. By $\dot{V}(x_{t})\mid_{\substack{ (*)  \\ u=u^{*}}}$, we denote the time derivative of the functional $V(x_{t})$ along the trajectories of system (*), when the control law is $u^{*}$.  

\section{Preliminaries and Problem Statement}
Consider time-delay systems of the form
\begin{equation}
\begin{split}
&\dot{x}(t)=Ax(t)+Bx(t-h)+Du(t),\\
&\varphi \in \mathrm{PC}([-h, 0],\mathbb{R}^{n})
\end{split}
\label{sistema1}
\end{equation}
where the matrices $A,B \in \mathbb{R}^{n \times n}$, $D\in \mathbb{R}^{n \times r}$ are constant, the state $x(t)$ is in $\mathbb{D}$, the space of solutions which contains the trivial one, and the control vector $u(t)$ belongs to $\mathbb{R}^{r}$, $r \leq n$.\\
Let the following quadratic performance index be given:
\begin{equation}
J=\int_{0}^{\infty} g(x_{t},u(t)) dt=\int\limits_0^\infty \left( x^{T}(t)Qx(t) + u^{T}(t)Ru(t)
\right) dt,
\label{desempenioJ}
\end{equation}
with $Q \in \mathbb{R}^{n \times n}, R \in \mathbb{R}^{r \times r}$, $Q=Q^{T} > 0,R=R^{T}>0$.\\
The optimal control problem consists in the synthesis of the optimal control $u^{*}(t)$ that minimizes the quadratic performance index (\ref{desempenioJ}) subject to \eqref{sistema1}. 

Admissible controls for this problem satisfy:
\begin{enumerate}
\item $u=u(x_{t})$, in other words, the control is a function of the system state, $x_{t}.$
\item The functional $u(x_{t})$\ is such that solutions to (\ref{sistema1}) exist and are unique for $t\geq 0$\ and for all initial conditions $\varphi $.
\item The trivial solution of (\ref{sistema1}) in closed-loop with control law $u=u(x_{t})$\
is asymptotically stable.
\item For $u=u(x_{t})$ and all initial conditions $\varphi $ the performance
index has a finite value. 
\end{enumerate}
 System (\ref{sistema1}) in closed-loop with an admissible control of the form
\begin{equation}
u_{L}(x_t)=\Gamma_{0}x(t)+\int\limits_{-h}^{0} \Gamma_{1}(\theta)x(t+\theta) d\theta,
\label{u_gammas}
\end{equation}
is an exponentially stable system given by
\begin{equation}
\dot{x}(t)= A_{0}x(t)+ A_{1}x(t-h)+\int\limits_{-h}^{0}  G(\theta)x(t+\theta) d\theta, t\geq 0.
\label{sistema_1_lazo_cerrado_u_gammas}
\end{equation}
where  $x(t) \in \mathbb{D}$; $A_{0}=A+D\Gamma_{0},A_{1}=B,G(\theta)=D\Gamma_{1}(\theta) \in \mathbb{R}^{n \times n}$.\\
\\The sufficient conditions showcasing the Bellman equation for an optimal control for time-delay system are given next:
\begin{thm}\rm{(Ross \cite{Ross_1969}).}
	\textit{If there exists an admissible control $u^{*}=u^{*}(x_{t})$ and a scalar continuous non negative function $V(x_{t})$, $V=0$ for all $x_{t}=0$ such that
	\begin{equation}
	\dot{V}(x_{t})\mid_{\substack{ (\ref{sistema1})  \\ u=u^{*}}} +
	g(x_{t},u^{*}(x_{t}))=0, \forall t\geq 0
	\label{inciso_a}
	\end{equation}
	\begin{equation}
	\dot{V}(x_{t}) \mid_{\substack{ (\ref{sistema1})  \\ u=u^{*}}} +
	g(x_{t},u^{*}(x_{t}))\\
	\leq \dot{V}(x_{t}) \mid_{\substack{ (\ref{sistema1})
			\\ u=u(t)}} + g(x_{t},u(t)), \forall t\geq 0
	\label{inciso_b}
	\end{equation}	
	for all admissible $u(t)$, then $u^{*}(t)$ is an optimal control. Furthermore $V(\varphi)=J(\varphi,u^{*})$ is the optimal value of  the performance index $J$.}
	\label{forma_Vx_t_ross}
\end{thm}
The functional $V(x_{t})$ is called the Bellman functional, which is used to provide the necessary and sufficient conditions of optimality for time-delay systems.

The necessary conditions for an optimal control for time-delay systems are given in the following proposition.

\begin{prop} \rm{(Ross \cite{Ross_1969}).}
	\textit{If $u_{L}=u_{L}(x_{t})$, $\forall t\geq 0$, is an admissible linear control, $\varphi$ is an
	initial condition function on $[-h,0]$,  then the function
	\begin{equation}\label{Vpositiva}
	V(\varphi) = J(\varphi, u_{L})= \int_{0}^{\infty}
	(x^{T}(t)Qx(t)+u_{L}(t)^{T}Ru_{L}(t)) dt,
	\end{equation}%
	can be expressed as
	\begin{equation}
	V(\varphi) = \varphi^{T}(0) \Pi_{0} \varphi(0) +2 \varphi^{T}(0) \int_{-h}^{0}
	\Pi_{1}(\theta) \varphi (\theta) d \theta
	+ \int_{-h}^{0} \int_{-h}^{0} \varphi^{T}(\xi) \Pi_{2} (\xi,\theta)
	\varphi(\theta) d \xi d \theta,
	\label{funcional_v_phi_11}
	\end{equation}
	where
	\begin{itemize}
		\item[i)] $\Pi_{0}>0$ is a symmetric positive matrix.	
		\item[ii)] $\Pi_{1}(\theta)$ is defined on $[-h,0]$.	
		\item[iii)] $\Pi_{2}(\xi,\theta)$ is defined on $%
		\xi, \theta \in [-h,0]$,\newline
		$\Pi_{2}^{T}(\xi, \theta)=\Pi_{2}(\theta,\xi)$.
	\end{itemize}
}
\label{estructura_v(phi)}
\end{prop}

The necessary and sufficient conditions for and optimal control for time-delay systems are given in the seminal result reminded bellow:
\begin{thm}\rm{(Ross \cite{Ross_1969}).}
	\textit{A linear control law
	\begin{equation}
	u^{*}(t)=-R^{-1}D^{T} \Pi_{0} x(t)
	-R^{-1}D^{T}\int_{-h}^{0}
	\Pi_{1}(\theta)x(t+ \theta) d\theta, \quad t\geq 0
	\label{ley_control_u0}
	\end{equation}
	provides the global minimum of the performance index \eqref{desempenioJ} for the dynamical system \eqref{sistema1} if:}
	\begin{itemize}
	    \item[a)] \textit{$u^{*}(x_{t})$ is a stabilizing control law (since $u^{*}$ is linear, stability and admissibility are equivalent)}
		\item[b)] \textit{$\Pi_{0}$ is a symmetric positive definite matrix which, together with the $n \times n$ array $\Pi_{1}(\theta)$ of functions defined on $[-h,0]$, and the $n\times n$ array, $\Pi_{2}(\xi,\theta)$ of functions in two variables having domain $\xi, \theta \in [-h,0]$, satisfies the relations:}
	\end{itemize}
	\begin{itemize}
		\item[1)] $
		A^{T}\Pi _{0}+\Pi _{0}A-\Pi _{0}DR^{-1}D^{T}\Pi _{0}+\Pi _{1}^{T}(0)$
		$+\Pi_{1}(0)+Q=0,$
		
		\item[2)] $\frac{d\Pi _{1}(\theta )}{d\theta }=(A^{T}-\Pi _{0}DR^{-1}D^{T})\Pi
		_{1}(\theta )+\Pi _{2}(0,\theta ),$
		$-h\leq \theta \leq 0,$
		
		\item[3)] $\frac{\partial \Pi _{2}(\xi ,\theta )}{\partial \xi }+\frac{\partial \Pi
			_{2}(\xi ,\theta )}{\partial \theta }=-\Pi _{1}^{T}(\xi )DR^{-h}D^{T}\Pi
		_{1}(\theta ),$
		$-h\leq \xi \leq 0,$ $-h\leq \theta \leq 0,$
		
		\item[4)] $\Pi _{1}(-h)=\Pi _{0}B,$
		
		\item[5)] $\Pi _{2}(-h,\theta )=B^{T}\Pi _{1}(\theta ),$ $-h\leq \theta \leq 0.$
		\begin{equation}
		\label{SetPartialEquations}
		\end{equation}
	\end{itemize}
	\textit{Furthermore, under these conditions, the representation of \eqref{desempenioJ} in
	terms of the initial function is}
	\begin{equation}
	\begin{split}
	&J(\varphi,u^{*})=\varphi^{T}(0) \Pi_{0} \varphi(0) + 2 \varphi^{T}(0)
	\int_{-h}^{0} \Pi_{1}(\theta) \varphi (\theta) d \theta\\ 
		&+ \int_{-h}^{0}	\int_{-h}^{0} \varphi^{T}(\xi) \Pi_{2} (\xi,\theta) \varphi(\theta) d \xi d
	\theta.
	\label{indice_desempenio_J_bajo_control_u0}
	\end{split}
	\end{equation} 
	\label{teorema_condiciones_necesarias y suficientes}
\end{thm} 
\section{Proof of some properties of the matrices in the Bellman functional}
In our contribution currently under revision \cite{Automatica2020}, we constructed the matrices $\Pi_{0}$, $\Pi_{1}(\theta)$ and $\Pi_{2}(\xi,\theta)$ that define the Bellman functional in Proposition \ref{estructura_v(phi)}. The obtained expressions in terms of the closed-loop system fundamental matrix $K(t)$ are as follows: 
\begin{equation}
    \begin{split}
        \Pi _{0} =&\int_{0}^{\infty }K^{T}(t)M_{1}K(t)dt+2\int_{-h}^{0}\left(
        \int_{0}^{\infty }K^{T}(t)M_{2}(\theta )K(t+\theta )dt\right) d\theta  \\
        &+\int_{-h}^{0}\int_{-h}^{0}\left( \int_{0}^{\infty }K^{T}(t+\theta
        _{1})M_{3}(\theta _{1},\theta _{2})K(t+\theta _{2})dt\right) d\theta
        _{2}d\theta _{1},    
    \end{split}
    \label{Pi_0_Ks_new}
\end{equation}
\begin{equation}
    \begin{split}
        \Pi _{1}(\theta ) =&\left( \int_{0}^{\infty }K^{T}(t)M_{1}K(t-\theta -h)dt%
        \right) A_{1}+\int_{-h}^{\theta }\left( \int_{0}^{\infty
        }K^{T}(t)M_{1}K(t-\theta +\xi )dt\right) G(\xi )d\xi  \\
        &+\int_{-h}^{0}\left( \int_{0}^{\infty }K^{T}(t)M_{2}(\theta
        _{2})K(t+\theta _{2}-\theta -h)dt\right) d\theta_{2} A_{1}
        +\int_{-h}^{0}\left( \int_{0}^{\infty }K^{T}(t+\theta
        _{2})M_{2}^{T}(\theta _{2})K(t-\theta -h)dt\right) d\theta _{2}A_{1}
        \\
        &+\int_{-h}^{0}\int_{-h}^{\theta }\left( \int_{0}^{\infty
        }K^{T}(t)M_{2}(\theta _{2})K(t+\theta _{2}-\theta +\xi )dt\right)
        G(\xi )d\xi d\theta _{2}\\
        &+\int_{-h}^{0}\int_{-h}^{\theta }\left(
        \int_{0}^{\infty }K^{T}(t+\theta _{2})M_{2}^{T}(\theta _{2})K(t-\theta
        +\xi )dt\right) G(\xi )d\xi d\theta _{2} \\
        &+\int_{-h}^{0}\int_{-h}^{0}\left( \int_{0}^{\infty }K^{T}(t+\theta
        _{1})M_{3}(\theta _{1},\theta _{2})K(t+\theta _{2}-\theta -h)dt\right) 
        d\theta _{2}d\theta_{1} A_{1}\\
        &+\int_{-h}^{0}\int_{-h}^{0}\int_{-h}^{\theta }\left( \int_{0}^{\infty
        }K^{T}(t+\theta _{1})M_{3}(\theta _{1},\theta _{2})K(t+\theta _{2}-\theta
        +\xi )dt \right) G(\xi )d\xi d\theta _{2}d\theta _{1},
    \end{split}
    \label{Pi_1_Ks_new}
\end{equation}
\begin{equation}
    \begin{split}
        \Pi _{2}(\xi ,\theta ) =&A_{1}^{T}\left( \int_{0}^{\infty
        }K^{T}(t-\xi -h)M_{1}K(t-\theta -h)dt\right) A_{1}
        +A_{1}^{T}\int_{-h}^{\theta
        }\left( \int_{0}^{\infty }K^{T}(t-\xi -h)M_{1}K(t-\theta
        +\delta )dt\right) G(\delta )d\delta  \\
        &+\int_{-h}^{\xi }G^{T}(\delta )\left( \int_{0}^{\infty }K^{T}(t-\xi
        +\delta )M_{1}K(t-\theta -h)dt\right) d\delta A_{1}\\ &+\int_{-h}^{\xi
        }\int_{-h}^{\theta }G^{T}(\delta _{1})\left( \int_{0}^{\infty }K^{T}(t-\xi
        +\delta _{1})M_{1}K(t-\theta +\delta _{2})dt\right) G(\delta _{2})d\delta
        _{2}d\delta _{1} \\
        &+2A_{1}^{T} \int_{-h}^{0}\left( \int_{0}^{\infty }K^{T}(t-\xi
        -h)M_{2}(\theta _{2})K(t+\theta _{2}-\theta -h)dt\right) d\theta _{2} A_{1}\\
        &+2A_{1}^{T}\int_{-h}^{0}\int_{-h}^{\theta }\left(
        \int_{0}^{\infty }K^{T}(t-\xi -h)M_{2}(\theta _{2})K(t+\theta _{2}-\theta
        +\delta )dt\right) G(\delta )d\delta d\theta _{2} \\
        &+2\int_{-h}^{0}\int_{-h}^{\xi }G^{T}(\delta )\left( \int_{0}^{\infty
        }K^{T}(t-\xi +\delta )M_{2}(\theta _{2})K(t+\theta _{2}-\theta -h)dt\right)
        d\delta d\theta _{2} A_{1} \\
        &+2\int_{-h}^{0}\int_{-h}^{\xi }\int_{-h}^{\theta }G^{T}(\delta _{1})\left(
        \int_{0}^{\infty }K^{T}(t-\xi +\delta _{1})M_{2}(\theta _{2})K(t+\theta
        _{2}-\theta +\delta _{2})dt\right) G(\delta _{2})d\delta _{2}d\delta
        _{1}d\theta _{2} \\
        &+A_{1}^{T}\int_{-h}^{0}\int_{-h}^{0}\left( \int_{0}^{\infty
        }K^{T}(t+\theta _{1}-\xi -h)M_{3}(\theta _{1},\theta _{2})K(t+\theta
        _{2}-\theta -h)dt\right) d\theta _{2}d\theta _{1} A_{1}\\
        &+A_{1}^{T}\int_{-h}^{0}\int_{-h}^{0}\int_{-h}^{\theta }\left(
        \int_{0}^{\infty }K^{T}(t+\theta _{1}-\xi -h)M_{3}(\theta _{1},\theta
        _{2})K(t+\theta _{2}-\theta +\delta )dt\right) G(\delta )d\delta d\theta
        _{2}d\theta _{1} \\
        &+\int_{-h}^{0}\int_{-h}^{0}\int_{-h}^{\xi }G^{T}(\delta )\left(
        \int_{0}^{\infty }K^{T}(t+\theta _{1}-\xi +\delta )M_{3}(\theta _{1},\theta
        _{2})K(t+\theta _{2}-\theta -h)dt\right) d\delta d\theta
        _{2}d\theta _{1} A_{1}\\
        &+\int_{-h}^{0}\int_{-h}^{0}\int_{-h}^{\xi }\int_{-h}^{\theta }G^{T}(\delta
        _{1})\left( \int_{0}^{\infty }K^{T}(t+\theta _{1}-\xi +\delta
        _{1})M_{3}(\theta _{1},\theta _{2})K(t+\theta _{2}-\theta +\delta _{2})dt%
        \right) G(\delta _{2})d\delta _{2}d\delta _{1}d\theta _{2}d\theta _{1},
    \end{split}
    \label{P2_1_Ks_new}
\end{equation}%
where    
	\begin{equation}
	    M_{1}= Q + \Gamma_{0}^{T}R\Gamma_{0},
	    \label{M1}
	\end{equation}
	\begin{equation}
	    M_{2}(\theta)=\Gamma_{0}^{T}R\Gamma_{1}(\theta), \quad \theta \in [-h,0],
	    \label{M2}
	\end{equation}
	\begin{equation}
	    M_{3}(\theta_{1},\theta_{2})=\Gamma_{1}^{T}(\theta_{1})R\Gamma_{1}(\theta_{2}), \quad \theta_{1}, \theta_{2} \in [-h,0].
	    \label{M3}
	\end{equation}\\
The above expressions fully justify the assumed form of the Bellman functional in Proposition 1. However, the properties i), ii) and iii) in Proposition \ref{estructura_v(phi)} \rm{(Ross \cite{Ross_1969})} of the matrices  $\Pi_{0}$, $\Pi_{1}(\theta)$ and $\Pi_{2}(\xi,\theta)$ are assumed to be true without proof. Below, we use expressions \eqref{Pi_0_Ks_new}- \eqref{M3} to formally prove that these properties are indeed satisfied.

\begin{prop}
 
The matrix $\Pi_{0}$ in Proposition \ref{estructura_v(phi)} is a symmetric positive matrix.
    \label{simetria_Pi_0}
\end{prop}
\textbf{\textit{Proof.}}
 Substituting \eqref{M1}-\eqref{M3} into  \eqref{Pi_0_Ks_new}, we get

\begin{equation*}
    \begin{split}
        \Pi _{0}=&
        \int_{0}^{\infty }\left( K^{T}(t)\left( Q+\Gamma_{0}^{T}R\Gamma_{0}\right)
        K(t)+2\int_{-h}^{0}K^{T}(t)\left( \Gamma_{0}^{T}R\Gamma_{1}(\theta )\right) K(t+\theta
        )d\theta \right. \\
        &+\left. \int_{-h}^{0}\int_{-h}^{0}K^{T}(t+\theta _{1})\left(
        \Gamma_{1}^{T}(\theta _{1})R\Gamma_{1}(\theta _{2})\right) K(t+\theta _{2})d\theta
        _{1}d\theta _{2}\right) dt.    
    \end{split}
\end{equation*}
As $\varphi^{T}(0) \Pi_{0} \varphi(0)$ is a scalar, the matrix \eqref{Pi_0_Ks_new} rewrites as:
\begin{equation*}
    \begin{split}
        \Pi _{0}=&
        \int_{0}^{\infty }\left( K^{T}(t)\left( Q+\Gamma_{0}^{T}R\Gamma_{0}\right)
        K(t)+\int_{-h}^{0}K^{T}(t)\Gamma_{0}^{T}R\Gamma_{1}(\theta )K(t+\theta )d\theta \right.
        +\int_{-h}^{0}K^{T}(t+\theta )\Gamma_{1}^{T}(\theta )R\Gamma_{0}K(t)d\theta \\
        &+\left. \left( \int_{-h}^{0}K^{T}(t+\theta )\Gamma_{1}^{T}(\theta )d\theta %
        \right) R\left( \int_{-h}^{0}\Gamma_{1}(\theta )K(t+\theta )d\theta \right) %
        \right) dt,
    \end{split}
\end{equation*}
hence $\Pi_{0}$ can be rewritten as the quadratic form 
\begin{equation}
\Pi _{0}%
\begin{array}{c}
=%
\end{array}%
\int_{0}^{\infty }\left[ 
\begin{array}{c}
K(t) \\ 
\int_{-h}^{0}\Gamma_{1}(\theta )K(t+\theta )d\theta 
\end{array}%
\right] ^{T} \left[ 
\begin{array}{cc}
Q+\Gamma_{0}^{T}R\Gamma_{0} & \Gamma_{0}^{T}R \\ 
R\Gamma_{0} & R%
\end{array}%
\right] \left[ 
\begin{array}{c}
K(t) \\ 
\int_{-h}^{0}\Gamma_{1}(\theta )K(t+\theta )d\theta%
\end{array}%
\right] dt,  \label{simetrica_Pi_0}
\end{equation}
As $R^{T}=R>0$ and $Q^{T}=Q>0$, and 
\begin{equation*}
\left( Q+\Gamma_{0}^{T}R\Gamma_{0}\right) -\left( \Gamma_{0}^{T}R\right) \left(
R^{-1}\right) \left( R\Gamma_{0}\right)>0,
\end{equation*}
 Schur complement lemma implies that 
\begin{equation}
\left[ 
\begin{array}{cc}
Q+\Gamma_{0}^{T}R\Gamma_{0} & \Gamma_{0}^{T}R \\ 
R\Gamma_{0} & R%
\end{array}%
\right]>0,
\end{equation}
thus $\Pi_0$ is a symmetric positive matrix. $\blacksquare$

\begin{prop}
The matrix $\Pi_{1}(\theta)$ in Proposition \ref{estructura_v(phi)} is defined on $[-h,0]$.
\end{prop}
\textbf{\textit{Proof.}} Consider the expression \eqref{Pi_1_Ks_new} for $\Pi_{1}(\theta)$. There, $K(t)$ is the fundamental matrix  of the closed-loop system \eqref{sistema_1_lazo_cerrado_u_gammas} corresponding to an admissible control, thus $K(t)$ is exponentially stable, hence all integral summands of \eqref{Pi_1_Ks_new} are defined on $[-h,0]$. $\blacksquare$

\begin{prop}
The matrix $\Pi_{2}(\theta,\xi)$ in Proposition \ref{estructura_v(phi)} is such that $\Pi_{2}^{T}(\xi, \theta)=\Pi_{2}(\theta,\xi)$.
\end{prop}
\textbf{\textit{Proof.}} 
The result is obtained by showing that \eqref{P2_1_Ks_new} and its transpose are equal. Observe first that the matrices $M_{1}$ and  $M_{3}(\theta_{1},\theta_{2})$ defined in \eqref{M1} and \eqref{M3} are such that 
\begin{equation}
	    M^{T}_{1}= Q + \Gamma_{0}^{T}R\Gamma_{0}=M_{1},
	    \label{simetria_M1}
\end{equation}
and
\begin{equation}
	    M_{3}^{T}(\theta_{1},\theta_{2})=\Gamma_{1}^{T}(\theta_{2})R\Gamma_{1}(\theta_{1})=M_{3}(\theta_{2},\theta_{1}).
	    \label{transpuesta_M3}
\end{equation}
Now, the fact that $ \int_{-h}^{0}	\int_{-h}^{0} \varphi^{T}(\xi) \Pi_{2} (\xi,\theta) \varphi(\theta) d \xi d\theta $ is a scalar, allows to verify that \eqref{P2_1_Ks_new} can be rewritten as:
\begin{equation}
\begin{split}
\Pi _{2}(\xi ,\theta ) =&\int_{0}^{\infty }\left( A_{1}^{T}K^{T}(t-\xi
        -h)M_{1}K(t-\theta -h)A_{1}+\int_{-h}^{\theta }A_{1}^{T}K^{T}(t-\xi
        -h)M_{1}K(t-\theta +\delta )G(\delta )d\delta \right. \\
 &+\int_{-h}^{\xi }G^{T}(\delta )K^{T}(t-\xi +\delta )M_{1}K(t-\theta
-h)A_{1}d\delta \\
 &+\int_{-h}^{\xi }\int_{-h}^{\theta }G^{T}(\delta _{1})K^{T}(t-\xi +\delta
_{1})M_{1}K(t-\theta +\delta _{2})G(\delta _{2})d\delta _{2}d\delta _{1} \\
&+\int_{-h}^{0}A_{1}^{T}K^{T}(t-\xi -h)M_{2}(\theta _{2})K(t+\theta _{2}-\theta-h)A_{1}d\theta _{2} \\
&+\int_{-h}^{0}A_{1}^{T}K^{T}(t+\theta _{2}-\xi -h)M_{2}^{T}(\theta
_{2})K(t-\theta -h)A_{1}d\theta _{2} \\
&+\int_{-h}^{0}\int_{-h}^{\theta }A_{1}^{T}K^{T}(t-\xi -h)M_{2}(\theta
_{2})K(t+\theta _{2}-\theta +\delta )G(\delta )d\delta d\theta _{2} \\
&+\int_{-h}^{0}\int_{-h}^{\theta }A_{1}^{T}K^{T}(t+\theta _{2}-\xi
-h)M_{2}^{T}(\theta _{2})K(t-\theta +\delta )G(\delta )d\delta d\theta _{2}
\\
&+\int_{-h}^{0}\int_{-h}^{\theta }\int_{-h}^{\xi }G^{T}(\delta
_{2})K^{T}(t+\theta _{2}-\xi +\delta _{2})M_{2}^{T}(\theta _{2})K(t-\theta
+\delta _{1})G(\delta _{1})d\delta _{2}d\delta _{1}d\theta _{2} \\
&+\int_{-h}^{0}\int_{-h}^{\xi }G^{T}(\delta )K^{T}(t+\theta _{2}-\xi
+\delta )M_{2}^{T}(\theta _{2})K(t-\theta -h)A_{1}d\delta d\theta _{2} \\
&+\int_{-h}^{0}\int_{-h}^{\xi }G^{T}(\delta )K^{T}(t-\xi +\delta
)M_{2}(\theta _{2})K(t+\theta _{2}-\theta -h)A_{1}d\delta d\theta _{2} \\
&+\int_{-h}^{0}\int_{-h}^{\xi }\int_{-h}^{\theta }G^{T}(\delta
_{1})K^{T}(t-\xi +\delta _{1})M_{2}(\theta _{2})K(t+\theta _{2}-\theta
+\delta _{2})G(\delta _{2})d\delta _{2}d\delta _{1}d\theta _{2} \\
&+\int_{-h}^{0}\int_{-h}^{0}A_{1}^{T}K^{T}(t+\theta _{1}-\xi -h)M_{3}(\theta
_{1},\theta _{2})K(t+\theta _{2}-\theta -h)A_{1}d\theta _{1}d\theta _{2} \\
&+\int_{-h}^{0}\int_{-h}^{0}\int_{-h}^{\theta }A_{1}^{T}K^{T}(t+\theta _{1}-\xi
-h)M_{3}(\theta _{1},\theta _{2})K(t+\theta _{2}-\theta +\delta )G(\delta
)d\delta d\theta _{1}d\theta _{2} \\
&+\int_{-h}^{0}\int_{-h}^{0}\int_{-h}^{\xi }G^{T}(\delta )K^{T}(t+\theta
_{1}-\xi +\delta )M_{3}(\theta _{1},\theta _{2})K(t+\theta _{2}-\theta
-h)A_{1}d\delta d\theta _{1}d\theta _{2} \\
&+\left. \int_{-h}^{0}\int_{-h}^{0}\int_{-h}^{\xi }\int_{-h}^{\theta
}G^{T}(\delta _{1})K^{T}(t+\theta _{1}-\xi +\delta _{1})M_{3}(\theta
_{1},\theta _{2})K(t+\theta _{2}-\theta +\delta _{2})G(\delta _{2})d\delta
_{2}d\delta _{1}d\theta _{1}d\theta _{2}\right) dt.
\end{split}
 \label{Pi_2_extendido}
\end{equation}
The expressions \eqref{P2_1_Ks_new} and \eqref{Pi_2_extendido} are equivalent, the only difference is that in \eqref{Pi_2_extendido} there are no sums of similar terms.
\newpage
If we transpose equation \eqref{Pi_2_extendido}, use the equalities \eqref{simetria_M1}, \eqref{transpuesta_M3} and carry out some changes of variables into the integrals, we get
\begin{equation}
    \begin{split}
        \Pi _{2}^{T}(\xi ,\theta ) =&\int_{0}^{\infty }\left( A_{1}^{T}K^{T}(t-\theta
        -h)M_{1}K(t-\xi -h)A_{1}+\int_{-h}^{\xi }A_{1}^{T}K^{T}(t-\theta -h)M_{1}K(t-\xi
        +\delta )G(\delta )d\delta \right.  \\
        &+\int_{-h}^{\theta }G^{T}(\delta )K^{T}(t-\theta +\delta )M_{1}K(t-\xi
        -h)A_{1} d\delta \\
        &+\int_{-h}^{\theta }\int_{-h}^{\xi }G^{T}(\delta _{1})K^{T}(t-\theta
        +\delta _{1})M_{1}K(t-\xi +\delta _{2})G(\delta _{2})d\delta _{2}d\delta _{1} \\
        &+\int_{-h}^{0}A_{1}^{T}K^{T}(t-\theta -h)M_{2}(\theta _{2})K(t+\theta _{2}-\xi
        -h)A_{1}d\theta _{2} \\
        &+\int_{-h}^{0}A_{1}^{T}K^{T}(t+\theta _{2}-\theta -h)M_{2}^{T}(\theta
        _{2})K(t-\xi -h)A_{1}d\theta _{2} \\
        &+\int_{-h}^{0}\int_{-h}^{\xi }A_{1}^{T}K^{T}(t-\theta -h)M_{2}(\theta
        _{2})K(t+\theta _{2}-\xi +\delta )G(\delta )d\delta d\theta _{2}\\
        &+\int_{-h}^{0}\int_{-h}^{\xi }A_{1}^{T}K^{T}(t+\theta _{2}-\theta
        -h)M_{2}^{T}(\theta _{2})K(t-\xi +\delta )G(\delta )d\delta d\theta _{2}\\
        &+\int_{-h}^{0}\int_{-h}^{\xi }\int_{-h}^{\theta }G^{T}(\delta
        _{2})K^{T}(t+\theta _{2}-\theta +\delta _{2})M_{2}^{T}(\theta _{2})K(t-\xi
        +\delta _{1})G(\delta _{1})d\delta _{2}d\delta _{1}d\theta _{2}\\
        &+\int_{-h}^{0}\int_{-h}^{\theta }G^{T}(\delta )K^{T}(t+\theta _{2}-\theta
        +\delta )M_{2}^{T}(\theta _{2})K(t-\xi -h)A_{1}d\delta d\theta _{2}\\
        &+\int_{-h}^{0}\int_{-h}^{\theta }G^{T}(\delta )K^{T}(t-\theta +\delta)M_{2}(\theta _{2})K(t+\theta _{2}-\xi -h)A_{1}d\delta d\theta _{2}\\
        &+\int_{-h}^{0}\int_{-h}^{\theta }\int_{-h}^{\xi }G^{T}(\delta
        _{1})K^{T}(t-\theta +\delta _{1})M_{2}(\theta _{2})K(t+\theta _{2}-\xi
        +\delta _{2})G(\delta _{2})d\delta _{2}d\delta _{1}d\theta _{2}\\
        &+\int_{-h}^{0}\int_{-h}^{0}A_{1}^{T}K^{T}(t+\theta _{1}-\theta -h)M_{3}(\theta
        _{1},\theta _{2})K(t+\theta _{2}-\xi -h)A_{1}d\theta _{1}d\theta _{2}\\
        &+\int_{-h}^{0}\int_{-h}^{0}\int_{-h}^{\xi }A_{1}^{T}K^{T}(t+\theta _{1}-\theta
        -h)M_{3}(\theta _{1},\theta _{2})K(t+\theta _{2}-\xi +\delta )G(\delta
        )d\delta d\theta _{1}d\theta _{2}\\
        &+\int_{-h}^{0}\int_{-h}^{0}\int_{-h}^{\theta }G^{T}(\delta )K^{T}(t+\theta
        _{1}-\theta +\delta )M_{3}(\theta _{1},\theta _{2})K(t+\theta _{2}-\xi
        -h)A_{1}d\delta d\theta _{1}d\theta _{2}\\
        &+\left. \int_{-h}^{0}\int_{-h}^{0}\int_{-h}^{\theta }\int_{-h}^{\xi
        }G^{T}(\delta _{1})K^{T}(t+\theta _{1}-\theta +\delta _{1})M_{3}(\theta
        _{1},\theta _{2})K(t+\theta _{2}-\xi +\delta _{2})G(\delta _{2})d\delta
        _{2}d\delta _{1}d\theta _{1}d\theta _{2}\right) dt.
    \end{split}
    \label{Pi_2_Trans}
\end{equation}

Finally, if $\xi$ and $\theta$ are interchanged in \eqref{Pi_2_Trans} then \eqref{Pi_2_extendido} is obtained, and the result follows. $\blacksquare$
\section{On the form of the optimal control for time-delay systems}
A natural query in the problem formulation is why the admissible controls are restricted to the form \eqref{u_gammas}. Indeed, because of the the delayed nature of the problem, it may seem natural to include as well a feedback term of the delayed state $x(t-h)$ in the control law.\\  
We remind next the procedure for determining optimal control suggested in \cite{Ross_1969}, page 615, line 12:\newline
\textit{(a) choose a particular form of }$u^{*}(x_{t});$\textit{\newline
(b) from that choice, express }$J(\varphi ,u^{*}(x_{t}))$ \textit{as an
explicit functional of the initial state, i.e., find }$V$\textit{\ such that
}$V(\varphi )=J(\varphi ,u^{*}(x_{t}))$\textit{;\newline
(c) use the equations \eqref{inciso_a}, \eqref{inciso_b}} \textit{as constraints on the
parameters of the assumed form of }$u^{*}(x_{t}).$\newline
At step (a) of the above algorithm, we choose the admissible control
\begin{equation}
u_{L}(x_{t})=\Gamma_{0}x(t)+\int_{-h}^{0}\Gamma_{1}(\theta )x(t+\theta )d\theta +\Gamma_{2}x(t-h).
\label{u_t}
\end{equation}
with matrix parameters $\Gamma_{0},\Gamma_{1}(\theta)$ and $\Gamma_{2},$ of appropriate dimensions. We need to find the matrix parameters that satisfy the sufficient conditions constraints (\ref{inciso_a}) and (%
\ref{inciso_b}) to conclude that $u^{*}(x_{t})\,$is optimal and $J(\varphi
,u^{*}(x_{t}))$ is a global minimum.

Substituting  the expression \eqref{u_t} into \eqref{sistema1}, gives the closed-loop system 
\begin{equation}
\dot{x}(t)=\tilde{A}_{0}x(t)+\tilde{A}_{1}x(t-h)+\int_{-h}^{0}\tilde{A}%
_{2}(\theta )x(t+\theta )d\theta,  \label{lazo_cerrado}
\end{equation}
where
\begin{equation*}
        \tilde{A}_{0} =A+D\Gamma_{0}, \quad
        \tilde{A}_{1} =B+D\Gamma_{2}, \quad
        \tilde{A}_{2}(\theta ) =D\Gamma_{1}(\theta ).        
\end{equation*}
The Cauchy formula \cite{Kolmanovskii_1992} for the distributed time-delay system \eqref{lazo_cerrado} is 
\begin{equation}
x(t,\varphi)=\tilde{K}(t) \varphi(0) + \int_{-h}^{0} \tilde{K}(t-\theta-h)\tilde{A}_{1}\varphi(\theta) d
\theta+ \int_{-h}^{0} \int_{-h}^{\theta} \tilde{K}(t- \theta+ \xi)\tilde{A}_{2}(\xi) d\xi
\varphi(\theta) d\theta, \quad t\geq 0,
\label{solucion_cauchy}
\end{equation}
where $\tilde{K}(t)$ is the fundamental matrix of system \eqref{lazo_cerrado}. It has exponentially stable trivial solution because the control law \eqref{u_t} is  admissible.

Introducing the control law \eqref{u_t} into the performance index \eqref{desempenioJ} yields
\begin{equation}
    \begin{split}
        J =&\int_{0}^{\infty}\bigg(x^{T}(t)L_{1}x(t)+2x^{T}(t)L_{2}x(t-h)+2\int_{-h}^{0}x^{T}(t)L_{3}(\theta)x(t+\theta)d\theta+x^{T}(t-h)L_{4}x(t-h)\\
        &+2\int_{-h}^{0}x^{T}(t-h)L_{5}(\theta)x(t+\theta)d\theta+\int_{-h}^{0}\int_{-h}^{0}x^{T}(t+\theta_{1})L_{6}(\theta_{1},\theta_{2})x(t+\theta_{2})d\theta_{2}d\theta_{1}\bigg) dt, 
    \end{split}
    \label{J_M_1_6}
\end{equation}
where
\begin{itemize}
	\setlength{\parskip}{0pt}
	\setlength{\itemsep}{0pt plus 1pt}
	\item $L_{1} =Q+\Gamma_{0}^{T}R\Gamma_{0}$,
	\item $L_{2}=\Gamma_{0}^{T}R\Gamma_{2}$,
	\item $L_{3}(\theta ) =\Gamma_{0}^{T}R\Gamma_{1}(\theta )$, $\theta \in [-h,0]$,
	\item $L_{4} =\Gamma_{2}^{T}R\Gamma_{2}$,
	\item $L_{5}(\theta ) =\Gamma_{2}^{T}R\Gamma_{1}(\theta )$, $\theta \in [-h,0]$,
	\item$L_{6}(\theta_{1},\theta_{2})=\Gamma_{1}^{T}(\theta_{1})R\Gamma_{1}(\theta _{2})$, $\theta_{1}, \theta_{2} \in [-h,0]$.
\end{itemize}
Substituting the Cauchy formula \eqref{solucion_cauchy} into \eqref{J_M_1_6}, we obtain

\begin{equation}
J=V(\varphi )=\varphi ^{T}(0)\tilde{\Pi} _{0}\varphi (0)+2\varphi
^{T}(0)\int_{-h}^{0}\tilde{\Pi} _{1}(\theta )\varphi (\theta )d\theta
+\int_{-h}^{0}\int_{-h}^{0}\varphi ^{T}(\xi )\tilde{\Pi} _{2}(\xi ,\theta )\varphi
(\theta )d\theta d\xi,
\label{J_igual_V}
\end{equation}
where
\begin{equation}
    \begin{split}
        \tilde{\Pi} _{0} =&\int_{0}^{\infty }\tilde{K}^{T}(t)L_{1}\tilde{K}(t)dt+2\int_{0}^{\infty
        }\tilde{K}^{T}(t)L_{2}\tilde{K}(t-h)dt+2\int_{-h}^{0}\left( \int_{0}^{\infty
        }\tilde{K}^{T}(t)L_{3}(\theta )\tilde{K}(t+\theta )dt\right) d\theta \\
        &+\int_{0}^{\infty }\tilde{K}^{T}(t-h)L_{4}\tilde{K}(t-h)dt+2\int_{-h}^{0}\left(
        \int_{0}^{\infty }\tilde{K}^{T}(t-h)L_{5}(\theta )\tilde{K}(t+\theta )dt\right) d\theta \\
        &+\int_{-h}^{0}\int_{-h}^{0}\left( \int_{0}^{\infty }\tilde{K}^{T}(t+\theta
        _{1})L_{6}(\theta _{1},\theta _{2})\tilde{K}(t+\theta _{2})dt\right) d\theta
        _{2}d\theta _{1},    
    \end{split}
    \label{Pi_0_Ks}
\end{equation}
\begin{equation}
    \begin{split}
        \tilde{\Pi} _{1}(\theta ) =&\left( \int_{0}^{\infty }\tilde{K}^{T}(t)L_{1}\tilde{K}(t-\theta -h)dt%
        \right) \tilde{A}_{1}+\int_{-h}^{\theta }\left(\int_{0}^{\infty
        }\tilde{K}^{T}(t)L_{1}\tilde{K}(t-\theta +\delta )dt\right) \tilde{A}_{2}(\delta )d\delta \\
        &+\left( \int_{0}^{\infty }\tilde{K}^{T}(t)L_{2}\tilde{K}(t-\theta -2h)dt\right) \tilde{A}%
        _{1}+\left( \int_{0}^{\infty }\tilde{K}^{T}(t-h)L_{2}^{T}\tilde{K}(t-\theta -h)dt\right)
        \tilde{A}_{1} \\
        &+\int_{-h}^{\theta }\left( \int_{0}^{\infty }\tilde{K}^{T}(t)L_{2}\tilde{K}(t-h-\theta
        +\delta )dt\right) \tilde{A}_{2}(\delta )d\delta +\int_{-h}^{\theta }\left(
        \int_{0}^{\infty }\tilde{K}^{T}(t-h)L_{2}^{T}\tilde{K}(t-\theta +\delta )dt\right) \tilde{A}_{2}(\delta
        )d\delta \\
        &+\int_{-h}^{0}\left( \int_{0}^{\infty }\tilde{K}^{T}(t)L_{3}(\theta
        _{2})\tilde{K}(t+\theta _{2}-\theta -h)dt\right) \tilde{A}_{1}d\theta
        _{2}+\int_{-h}^{0}\left( \int_{0}^{\infty }\tilde{K}^{T}(t+\theta
        _{2})L_{3}^{T}(\theta _{2})\tilde{K}(t-\theta -h)dt\right) \tilde{A}_{1}d\theta _{2}
        \\
        &+\int_{-h}^{0}\int_{-h}^{\theta }\left( \int_{0}^{\infty
        }\tilde{K}^{T}(t)L_{3}(\theta _{2})\tilde{K}(t+\theta _{2}-\theta +\delta )dt\right)
        \tilde{A}_{2}(\delta )d\delta d\theta _{2}\\
        &+\int_{-h}^{0}\int_{-h}^{\theta }\left(
        \int_{0}^{\infty }\tilde{K}^{T}(t+\theta _{2})L_{3}^{T}(\theta _{2})\tilde{K}(t-\theta
        +\delta )dt\right) \tilde{A}_{2}(\delta )d\delta d\theta _{2} \\
        &+\left( \int_{0}^{\infty }\tilde{K}^{T}(t-h)L_{4}\tilde{K}(t-\theta -2h)dt\right) \tilde{A}%
        _{1}+\int_{-h}^{\theta }\left( \int_{0}^{\infty }\tilde{K}^{T}(t-h)L_{4}\tilde{K}(t-h-\theta
        +\delta )dt\right) \tilde{A}_{2}(\delta )d\delta \\
        &+\int_{-h}^{0}\left( \left( \int_{0}^{\infty }\tilde{K}^{T}(t-h)L_{5}(\theta
        _{2})\tilde{K}(t+\theta _{2}-\theta -h)dt\right) \tilde{A}_{1}\right) d\theta
        _{2}\\
        &+\int_{-h}^{0}\left( \int_{0}^{\infty }\tilde{K}^{T}(t+\theta
        _{2})L_{5}^{T}(\theta _{2})\tilde{K}(t-\theta -2h)dt\right) \tilde{A}_{1}d\theta _{2}
        \\  
        &+\int_{-h}^{0}\int_{-h}^{\theta }\left( \int_{0}^{\infty
        }\tilde{K}^{T}(t-h)L_{5}(\theta _{2})\tilde{K}(t+\theta _{2}-\theta +\delta )dt\right)
        \tilde{A}_{2}(\delta )d\delta d\theta _{2}\\
        &+\int_{-h}^{0}\int_{-h}^{\theta }\left(
        \int_{0}^{\infty }\tilde{K}^{T}(t+\theta _{2})L_{5}^{T}(\theta _{2})\tilde{K}(t-h-\theta
        +\delta )dt\right) \tilde{A}_{2}(\delta )d\delta d\theta _{2} \\
        &+\int_{-h}^{0}\int_{-h}^{0}\left( \int_{0}^{\infty }\tilde{K}^{T}(t+\theta
        _{1})L_{6}(\theta _{1},\theta _{2})\tilde{K}(t+\theta _{2}-\theta -h)dt\right) 
        \tilde{A}_{1}d\theta _{2}d\theta
        _{1}\\
        &+\int_{-h}^{0}\int_{-h}^{0}\int_{-h}^{\theta }\left( \int_{0}^{\infty
        }\tilde{K}^{T}(t+\theta _{1})L_{6}(\theta _{1},\theta _{2})\tilde{K}(t+\theta _{2}-\theta
        +\delta )dt\right) \tilde{A}_{2}(\delta )d\delta d\theta _{2}d\theta _{1},
    \end{split}
    \label{Pi_1_Ks}
\end{equation}
and
\begin{equation*}
    \begin{split}
        \tilde{\Pi} _{2}(\xi ,\theta ) =&\tilde{A}_{1}^{T}\left( \int_{0}^{\infty
        }\tilde{K}^{T}(t-\xi -h)L_{1}\tilde{K}(t-\theta -h)dt\right) \tilde{A}_{1}+\int_{-h}^{\theta
        }\tilde{A}_{1}^{T}\left( \int_{0}^{\infty }\tilde{K}^{T}(t-\xi -h)L_{1}\tilde{K}(t-\theta
        +\delta )dt\right) \tilde{A}_{2}(\delta )d\delta \\
        &+\int_{-h}^{\xi }\tilde{A}_{2}^{T}(\delta )\left( \int_{0}^{\infty }\tilde{K}^{T}(t-\xi
        +\delta )L_{1}\tilde{K}(t-\theta -h)dt\right) \tilde{A}_{1}d\delta\\
        &+\int_{-h}^{\xi
        }\int_{-h}^{\theta }\tilde{A}_{2}^{T}(\delta _{1})\left( \int_{0}^{\infty }\tilde{K}^{T}(t-\xi
        +\delta _{1})L_{1}\tilde{K}(t-\theta +\delta _{2})dt\right) \tilde{A}_{2}(\delta _{2})d\delta
        _{2}d\delta _{1} \\
        &+2\tilde{A}_{1}^{T}\left( \int_{0}^{\infty }\tilde{K}^{T}(t-\xi -h)L_{2}\tilde{K}(t-\theta
        -2h)dt\right) \tilde{A}_{1}\\
        &+2\int_{-h}^{\theta }\tilde{A}_{1}^{T}\left(
        \int_{0}^{\infty }\tilde{K}^{T}(t-\xi -h)L_{2}\tilde{K}(t-h-\theta +\delta )dt\right)
        \tilde{A}_{2}(\delta )d\delta \\
        &+2\int_{-h}^{\xi }\tilde{A}_{2}^{T}(\delta )\left( \int_{0}^{\infty }\tilde{K}^{T}(t-\xi
        +\delta)L_{2}\tilde{K}(t-\theta -2h)dt\right) \tilde{A}_{1}d\delta\\
        &+2\int_{-h}^{\xi
        }\int_{-h}^{\theta }\tilde{A}_{2}^{T}(\delta _{1})\left( \int_{0}^{\infty }\tilde{K}^{T}(t-\xi
        +\delta _{1})L_{2}\tilde{K}(t-h-\theta +\delta _{2})dt\right) \tilde{A}_{2}(\delta _{2})d\delta
        _{2}d\delta _{1} \\
        \end{split}
\end{equation*}
\begin{equation}
    \begin{split}
         &+2\int_{-h}^{0}\tilde{A}_{1}^{T}\left( \int_{0}^{\infty }\tilde{K}^{T}(t-\xi
        -h)L_{3}(\theta _{2})\tilde{K}(t+\theta _{2}-\theta -h)dt\right) \tilde{A}%
        _{1}d\theta _{2}\\
        &+2\int_{-h}^{0}\int_{-h}^{\theta }\tilde{A}_{1}^{T}\left(
        \int_{0}^{\infty }\tilde{K}^{T}(t-\xi -h)L_{3}(\theta _{2})\tilde{K}(t+\theta _{2}-\theta
        +\delta )dt\right) \tilde{A}_{2}(\delta )d\delta d\theta _{2} \\
        &+2\int_{-h}^{0}\int_{-h}^{\xi }\tilde{A}_{2}^{T}(\delta )\left( \int_{0}^{\infty
        }\tilde{K}^{T}(t-\xi +\delta )L_{3}(\theta _{2})\tilde{K}(t+\theta _{2}-\theta -h)dt\right) 
        \tilde{A}_{1}d\delta d\theta _{2} \\
        &+2\int_{-h}^{0}\int_{-h}^{\xi }\int_{-h}^{\theta }\tilde{A}_{2}^{T}(\delta _{1})\left(
        \int_{0}^{\infty }\tilde{K}^{T}(t-\xi +\delta _{1})L_{3}(\theta _{2})\tilde{K}(t+\theta
        _{2}-\theta +\delta _{2})dt\right) \tilde{A}_{2}(\delta _{2})d\delta _{2}d\delta
        _{1}d\theta _{2} \\
        &+\tilde{A}_{1}^{T}\left( \int_{0}^{\infty }\tilde{K}^{T}(t-\xi -2h)L_{4}K(t-\theta
        -2h)dt\right) \tilde{A}_{1}\\
        &+\int_{-h}^{\theta }\tilde{A}_{1}^{T}\left(
        \int_{0}^{\infty }\tilde{K}^{T}(t-\xi -2h)L_{4}\tilde{K}(t-h-\theta +\delta )dt\right)
        \tilde{A}_{2}(\delta )d\delta \\
        &+\int_{-h}^{\xi }\tilde{A}_{2}^{T}(\delta )\left( \int_{0}^{\infty }\tilde{K}^{T}(t-h-\xi
        +\delta )L_{4}\tilde{K}(t-\theta -2h)dt\right) \tilde{A}_{1}d\delta \\
        &+\int_{-h}^{\xi }\int_{-h}^{\theta }\tilde{A}_{2}^{T}(\delta _{1})\left(
        \int_{0}^{\infty }\tilde{K}^{T}(t-h-\xi +\delta _{1})L_{4}\tilde{K}(t-h-\theta +\delta
        _{2})dt\right) \tilde{A}_{2}(\delta _{2})d\delta _{2}d\delta _{1}\\
        &+2\int_{-h}^{0}\tilde{A}_{1}^{T}\left( \int_{0}^{\infty }\tilde{K}^{T}(t-\xi
        -2h)L_{5}(\theta _{2})\tilde{K}(t+\theta _{2}-\theta -h)dt\right) \tilde{A}%
        _{1}d\theta _{2} \\
        &+2\int_{-h}^{0}\int_{-h}^{\theta }\tilde{A}_{1}^{T}\left( \int_{0}^{\infty
        }\tilde{K}^{T}(t-\xi -2h)L_{5}(\theta _{2})\tilde{K}(t+\theta _{2}-\theta +\delta )dt\right)
        \tilde{A}_{2}(\delta )d\delta d\theta _{2} \\
        &+2\int_{-h}^{0}\int_{-h}^{\xi }\tilde{A}_{2}^{T}(\delta )\left( \int_{0}^{\infty
        }\tilde{K}^{T}(t-h-\xi +\delta )L_{5}(\theta _{2})\tilde{K}(t+\theta _{2}-\theta -h)dt\right)
        \tilde{A}_{1}d\delta d\theta _{2} \\
        &+2\int_{-h}^{0}\int_{-h}^{\xi }\int_{-h}^{\theta }\tilde{A}_{2}^{T}(\delta _{1})\left(
        \int_{0}^{\infty }\tilde{K}^{T}(t-h-\xi +\delta _{1})L_{5}(\theta _{2})\tilde{K}(t+\theta
        _{2}-\theta +\delta _{2})dt\right) \tilde{A}_{2}(\delta _{2})d\delta _{2}d\delta
        _{1}d\theta _{2} \\
        &+\int_{-h}^{0}\int_{-h}^{0}\tilde{A}_{1}^{T}\left( \int_{0}^{\infty
        }\tilde{K}^{T}(t+\theta _{1}-\xi -h)L_{6}(\theta _{1},\theta _{2})\tilde{K}(t+\theta
        _{2}-\theta -h)dt\right) \tilde{A}_{1}d\theta _{2}d\theta _{1} \\
        &+\int_{-h}^{0}\int_{-h}^{0}\int_{-h}^{\theta }\tilde{A}_{1}^{T}\left(
        \int_{0}^{\infty }\tilde{K}^{T}(t+\theta _{1}-\xi -h)L_{6}(\theta _{1},\theta
        _{2})\tilde{K}(t+\theta _{2}-\theta +\delta )dt\right) \tilde{A}_{2}(\delta )d\delta d\theta
        _{2}d\theta _{1} \\
        &+\int_{-h}^{0}\int_{-h}^{0}\int_{-h}^{\xi }\tilde{A}_{2}^{T}(\delta )\left(
        \int_{0}^{\infty }\tilde{K}^{T}(t+\theta _{1}-\xi +\delta )L_{6}(\theta _{1},\theta
        _{2})\tilde{K}(t+\theta _{2}-\theta -h)dt\right) \tilde{A}_{1}d\delta d\theta
        _{2}d\theta _{1} \\
        &+\int_{-h}^{0}\int_{-h}^{0}\int_{-h}^{\xi }\int_{-h}^{\theta }\tilde{A}_{2}^{T}(\delta
        _{1})\left( \int_{0}^{\infty }\tilde{K}^{T}(t+\theta _{1}-\xi +\delta
        _{1})L_{6}(\theta _{1},\theta _{2})\tilde{K}(t+\theta _{2}-\theta +\delta _{2})dt%
        \right) \tilde{A}_{2}(\delta _{2})d\delta _{2}d\delta _{1}d\theta _{2}d\theta _{1}.
    \end{split}
    \label{P2_1_Ks}
\end{equation}%
At this point, we conclude that the form of the functional for the control \eqref{u_t} is also a tree term functional of the form \eqref{J_igual_V}, now with the above defined matrices. Its expression in terms of the state $x_{t}$ is
\begin{equation} 
V(x_{t})=x^{T}(t)\tilde{\Pi} _{0}x(t)+2x^{T}(t)\int_{-h}^{0}\tilde{\Pi} _{1}(\theta
)x(t+\theta )d\theta +\int_{-h}^{0}\int_{-h}^{0}x^{T}(t+\xi )\tilde{\Pi} _{2}(\xi
,\theta )x(t+\theta )d\theta d\xi ,\quad \forall t\geq 0.
\label{funcional_V_x_t}
\end{equation}%
Suppose that $V(x_{t})\geq 0$ and $x_{t}$ is a trajectory of system (\ref%
{sistema1}). Let us define
\begin{equation}
H(x_{t},u)=\left. \dot{V}(x_{t})\right\vert _{\substack{ (\ref{sistema1})  \\ %
u-admissible}}+x^{T}(t)Qx(t)+u^{T}(t)Ru(t),  \label{H_x_u}
\end{equation}%
where the time derivative of \eqref{funcional_V_x_t} along the trajectories of the system (\ref%
{sistema1}) is
\begin{equation}
    \begin{split}
        &\left. \dot{V}(x_{t})\right\vert _{\substack{ (\ref{sistema1})  \\ u-admissible}}
        =\left( Ax(t)+Bx(t-h)+Du(t)\right) ^{T}\tilde{\Pi} _{0}x(t)+x^{T}(t)\tilde{\Pi}
        _{0}\left( Ax(t)+Bx(t-h)+Du(t)\right) \\
        &+\left( Ax(t)+Bx(t-h)+Du(t)\right) ^{T}\int_{-h}^{0}\tilde{\Pi}
        _{1}(\theta )x(t+\theta )d\theta +x^{T}(t)\int_{-h}^{0}\tilde{\Pi} _{1}(\theta )%
        \frac{\partial }{\partial \theta }x(t+\theta )d\theta \\
        &+\int_{-h}^{0}\frac{\partial }{\partial \theta }x^{T}(t+\theta )\tilde{\Pi}
        _{1}^{T}(\theta )d\theta x(t)+\int_{-h}^{0}x^{T}(t+\theta )\tilde{\Pi}
        _{1}^{T}(\theta )d\theta \left( Ax(t)+Bx(t-h)+Du(t)\right) \\
        &+\int_{-h}^{0}\int_{-h}^{0}\frac{\partial }{\partial \xi }\left(
        x^{T}(t+\xi )\right) \tilde{\Pi} _{2}(\xi ,\theta )x(t+\theta )d\theta d\xi
        +\int_{-h}^{0}\int_{-h}^{0}x^{T}(t+\xi )\tilde{\Pi} _{2}(\xi ,\theta )\frac{\partial 
        }{\partial \theta }\left( x(t+\theta )\right) d\theta d\xi.
    \end{split}
    \label{V_dot_x_t}
\end{equation}
Substituting \eqref{V_dot_x_t} into \eqref{H_x_u}, implies that
\begin{equation}
\begin{split}
H(x_{t},u) =&x^{T}(t)\left( A^{T}\tilde{\Pi} _{0}+\tilde{\Pi} _{0}A+Q\right)
x(t)+2x^{T}(t)\tilde{\Pi} _{0}Bx(t-h)+2x^{T}(t)\tilde{\Pi} _{0}Du(t) \\
&+2x^{T}(t)A^{T}\int_{-h}^{0}\tilde{\Pi} _{1}(\theta )x(t+\theta )d\theta
+2x^{T}(t-h)B^{T}\int_{-h}^{0}\tilde{\Pi} _{1}(\theta )x(t+\theta )d\theta \\
&+2u^{T}(t)D^{T}\int_{-h}^{0}\tilde{\Pi} _{1}(\theta )x(t+\theta )d\theta
+x^{T}(t)\int_{-h}^{0}\tilde{\Pi} _{1}(\theta )\frac{\partial }{\partial \theta }%
x(t+\theta )d\theta  \\
&+\int_{-h}^{0}\frac{\partial }{\partial \theta }x^{T}(t+\theta )\tilde{\Pi}
_{1}^{T}(\theta )d\theta x(t)+\int_{-h}^{0}\int_{-h}^{0}\frac{\partial }{%
\partial \xi }[x^{T}(t+\xi )]\tilde{\Pi} _{2}(\xi ,\theta )x(t+\theta )d\theta d\xi 
\\
&+\int_{-h}^{0}\int_{-h}^{0}x^{T}(t+\xi )\tilde{\Pi} _{2}(\xi ,\theta )\frac{%
\partial }{\partial \theta }[x(t+\theta )]d\theta d\xi +u^{T}(t)Ru(t).
\end{split}
\label{H_x_t_u_1}
\end{equation}%
By the fundamental theorem of calculus of variations \cite{kirk2004optimal}
\begin{equation*}
\min_{u-admissible}H(x_{t},u)=H(x_{t}^{\ast },u^{\ast }),
\end{equation*}%
\begin{equation*}
\frac{\partial }{\partial u}H(x_{t},u)=2D^{T}\tilde{\Pi} _{0}x(t)+2D^{T}\int_{-h}^{0}\tilde{\Pi} _{1}(\theta )x(t+\theta )d\theta+2Ru(t)=0
\end{equation*}%
then, we have that 
\begin{equation}
u^{*}(t)=-R^{-1}D^{T}\tilde{\Pi} _{0}x(t)-R^{-1}D^{T}\int_{-h}^{0}\tilde{\Pi} _{1}(\theta
)x(t+\theta )d\theta.  \label{u_0_t_optima}
\end{equation}
Moreover, as
\begin{equation*}
\frac{\partial ^{2}}{\partial u^{2}}H(x_{t},u)=2R>0.
\end{equation*}%
we conclude that $u^{*}(t)$ is a local minimum of \eqref{H_x_u}.
Therefore the optimal gain for the $x(t-h)$ term in the admissible control law \eqref{u_t} is $\Gamma_2=0$, consequently, $L_{2}=0$, $ L_{4} =0,$ and $ L_{5}(\theta ) =0$. Furthermore, the expressions for $\tilde{\Pi} _{0}$, $\tilde{\Pi} _{1}(\theta )$  $\tilde{\Pi} _{2}(\xi ,\theta )$ reduce to \eqref{Pi_0_Ks_new}, \eqref{Pi_1_Ks_new}, \eqref{P2_1_Ks_new},  respectively, with $L_{1}=M_{1}$, $L_{3}(\theta)=M_{2}(\theta),$ and $L_{6}(\theta_{1},\theta_{2})=M_{3}(\theta_{1},\theta_{2})$. \textbf{By using the sufficient conditions \eqref{inciso_a} and \eqref{inciso_b}, it is possible to prove that the control \eqref{u_0_t_optima} is a global optimal control. }
\section{Concluding remarks}
In this note, we have proved some interesting complementary results concerning some properties of the Bellman functional, and on the form of the optimal control. It appears that establishing formally that the functional has a given structure is a crucial step in the solution of each of these problems. 
\bibliographystyle{unsrt}  
\bibliography{references}  
\end{document}